\documentclass[12pt]{amsart}

\def\@typesizes{%
       \or{5}{6.5}\or{6}{7.5}\or{7}{8.5}\or{8}{11}\or{9}{12}%
       \or{10}{13}
       \or{\@xipt}{14}\or{\@xiipt}{15}\or{\@xivpt}{18}%
       \or{\@xviipt}{20}\or{\@xxpt}{24}}

\usepackage{amsthm}
\usepackage{amsmath}
\usepackage{amssymb}
\usepackage{graphicx}
\usepackage{enumerate}
\allowdisplaybreaks
\numberwithin{equation}{section}
\numberwithin{figure}{section}

\theoremstyle{plain}
\newtheorem{theorem}{ Theorem}[section]
\newtheorem{proposition}[theorem]{ Proposition}
\newtheorem{lemma}[theorem]{ Lemma}
\newtheorem{corollary}[theorem]{ Corollary}
\newtheorem{example}[theorem]{ Example}
\newtheorem{remark}[theorem]{ Remark}
\newtheorem{definition}[theorem]{ Definition}
\newtheorem{conjecture}{ Conjecture}

\def\BET{\begin{theorem}}
\def\ENT{\end{theorem}}
\def\BEP{\begin{proposition}}
\def\ENP{\end{proposition}}
\def\BEL{\begin{lemma}}
\def\ENL{\end{lemma}}
\def\BEC{\begin{corollary}}
\def\ENC{\end{corollary}}
\def\BEE{\begin{example} \rm}
\def\ENE{\end{example}}
\def\BER{\begin{remark} \rm}
\def\ENR{\end{remark}}
\def\BED{\begin{definition} \rm}
\def\END{\end{definition}}
\def\BECJ{\begin{conjecture}}
\def\ENCJ{\end{conjecture}}

\def\bea{\begin{eqnarray}}
\def\eea{\end{eqnarray}}

\def\beas{\begin{eqnarray*}}
\def\eeas{\end{eqnarray*}}

\def\beq{\begin{equation}}
\def\eeq{\end{equation}}

\def\beal{\begin{align*}}

\def\eeal{ \end{align*} }

\def\roweq{\nonumber \\ &=& }
\def\rowleq{\nonumber \\  & \leq & }

\def\bfR{{\bf R}}

\def\bbC{{\mathbb C}}

\def\bbK{{\mathbb K}}

\def\bbN{{\mathbb N}}

\def\bbR{{\mathbb R}}
\def\bbS{{\mathbb S}}

\begin{document}

\title{Schauder bases and the decay rate of the heat equation}

\author{Jos\'e Bonet}
\author{Wolfgang Lusky}
\author{Jari Taskinen}



\begin{abstract}
We consider the classical Cauchy problem for the linear heat equation and integrable 
initial data in the Euclidean space $\bbR^N$. In the case $N=1$ we show that given a 
weighted $L^p$-space  $L_w^p(\bbR)$ with $1 \leq p < \infty$ and a fast growing weight 
$w$, there is  a Schauder basis $(e_n)_{n=1}^\infty$ in $L_ w^p(\bbR)$ with the 
following  property: given a positive 
integer $m $ there exists $n_m > 0$ such that, if the initial data $f$ belongs to the 
closed linear space of $e_n$ with $n \geq n_m$, then the decay rate of the solution 
of the heat equation is at least $t^{-m}$.  The result is also generalized to the case $N 
>1$ with a slightly weaker formulation. The proof is based on a construction 
of a Schauder basis of $L_w^p( \bbR^N)$,  which annihilates an infinite sequence of bounded 
functionals.

\end{abstract}

\maketitle

\section{Introduction.}
\label{sec1}

Given an integrable function $f \in L^1(\bbR^N)$ in the Euclidean space
$\bbR^N$, $N\in \bbN$, the unique solution
of the  classical Cauchy problem for the linear heat (or diffusion) equation
\bea
\partial_t u(x,t) &=& \Delta u(x,t) \ \ \mbox{for $x\in \bbR^N$, $t>0$} \label{1.2} \\
u(x,0) &=& f(x) \ \ \mbox{for $x\in \bbR^N$},\label{1.3}
\eea
has the decay rate $t^{-N/2}$ for large "times" $t$. This follows directly from the well-known  solution formula
\bea
u(x,t)= e^{t \Delta } f(x) := \frac{1}{(2\pi t)^{N/2}}
\int\limits_{\bbR^N} e^{- \frac1{4t} (x-y)^2 } f(y) dy , \label{1.4}
\eea
where we write $x^2 := |x|^2 = \sum_{j=1}^N x_j^2$ for vectors $x =(x_1, \ldots ,x_N) 
\in \bbR^N$ and $\Delta = \sum_{j=1}^N \partial_j^2  = \sum_{j=1}^N 
(\partial / \partial x_j)^2 $ for the Laplacian; if $N=1$ we denote $ e^{t \partial_x^2 } f$
instead of $ e^{t \Delta } f$. Indeed, \eqref{1.4} implies the bound 
\bea
\Vert u( \cdot,t ) \Vert_{p} := \Big( \int\limits_{\bbR^N} |u(x,t)|^p dx \Big)^{1/p} 
\leq C_p t^{-N/2}
\eea
for large $t$, for any $p \in [1, \infty)$ and also the same estimate for the sup-norm
$\Vert u(\cdot,t) \Vert_\infty$ with the usual definition. 

For general initial data $f \in L^1(\bbR^N)$, which is not necessarily positive, cancellation
phenomena may cause faster decay rates. For example in the case
$N=1$, if $f $ is such that $\int_{-\infty}^\infty f (x) dx =0$, then  a simple
argument shows that $e^{t \partial_x^2 } f$ decays at least with the speed $t^{-1}$; see  
Proposition  \ref{prop1.3} for an exact, more general formulation of this phenomenon. 

To describe our main result on decay rates we  fix a continuous weight function $w : 
\bbR^N \to \bbR^+ $ with symmetry  $w(x) := w(-x)$ for 
all $x \in \bbR^N$. We assume $w$ is fast growing  which means that  
\bea
\sup\limits_{x \in \bbR^N} \frac1{w(x)} (1 + |x|)^m < \infty \ \ \forall \, m \in \bbN .
\label{1.6}
\eea
Given $p \in [1, \infty)$ we denote by $
L_w^p(\bbR^N)$ the weighted $L^p$-space on $\bbR^N$ endowed with the norm
\bea
\Vert f \Vert_{p,w} := \Big( \int\limits_{\bbR^N} |f(x)|^p w(x) dx  \Big)^{1/p} .
\eea

Our main result, in addition to Theorem \ref{th1.4} on Schauder bases which annihilate 
linear functionals, reads as follows:

\BET
\label{th1.0}
Let $1 \leq p < \infty$ and let the  weight $w$ satisfy the conditions above. 

\smallskip

\noindent
$1^\circ$. Let $N=1$.  There exists a Schauder basis
$(e_n)_{n=1}^\infty$ of the Banach-space $L_w^p(\bbR)$  with the following
property: given  $m \in \bbN$ there exists $n_m \in \bbN$ such that any initial data $f$
\bea
f = \sum_{n=1}^{\infty} f_n e_n   \in L_w^p(\bbR)  ,  \label{1.6y}
\eea
with the property $f_n = 0$ for all $n=1,\ldots, n_m$, has the fast decay property
\bea
\Vert  e^{t \Delta} f \Vert_\infty \leq \frac{C}{t^{m}} \Vert f \Vert_{p,w} \label{1.12a}
\eea
for all $t\geq 1 $. 

\smallskip

\noindent $2^\circ$. If $N > 1$, then there exists a weight $\tilde w : \bbR^N \to \bbR^+$
satisfying the assumptions around \eqref{1.6} such that $ L^p_{w}(\bbR^N)  \subset L^1_{\tilde w}(\bbR^N)$ and such that 
the space $L_{\tilde w}^1(\bbR^N)$ has a Schauder basis $(e_n)_{n=1}^\infty$
with  the same property as in $1^\circ$ 
($\Vert f \Vert_{1, \tilde w}$ replacing $\Vert f \Vert_{p,w}$ in \eqref{1.12a}) .
\ENT

In other words, if initial data is included in the finite co-dimensional subspace
$G_m = \overline{{\rm sp} (e_n  :  n \geq n_m ) }$, then the corresponding 
solution decays at least at the speed $t^{-m}$; leaving out finitely many coordinates 
in the Banach-space of initial data makes the solution decay fast. 
The subspace $G_m$ thus has an explicit  description in terms of the Schauder basis, 
although in  general we are not able to determine the precise decay rate, if the  initial 
data is in the complement space. 

If $N > 1$ and $p=1$ and  the weight $w$ has a special symmetric form, then we may
still take $\tilde w = w$. See Section \ref{sec1.b} for details. 

\BER
a) We emphasize the following general aspect of our construction in the case $N=1$,
$p>1$:  to find the basis we split the space $L_w^p(\bbR) = L_w^{p,-}(\bbR) \oplus
L_w^{p,+}(\bbR) $, where the two subspaces consist of functions vanishing on the
positive or negative real line, respectively. Then, the basis in  Theorem \ref{th1.0}, 
is constructed as  small perturbations of any given Schauder bases
of $L_w^{p,\pm}(\bbR) $. Due this general nature of the 
result, we only  obtain  the existence of the desired basis, but not  explicit
information on the magnitude of the numbers $n_m$. See the end of  Section \ref{sec1.b}.

b) By classical arguments, the heat kernel in \eqref{1.4} can be expanded as the series
\bea
e^{\frac1{4t}(x-y)^2}  = \sum_{n \in  \bbN_0} \frac1{t^{n/2}} H_n(x) y^n 
\eea 
where $H_n$ are suitably normalized Hermite functions. If $m \in \bbN$,
one can write a given $f$, say, belonging to $L_w^2(\bbR)$ with $w(x) = e^{-x^2/2}$,  as 
\beas
f = \sum_{n=1}^m f_n H_n  + g ,
\eeas
where the coefficients $f_n$ are chosen such that  $\int_\bbR y^n g(y) dy = 0$
for $n=1, \ldots, m$.
Then, the solution with initial data  $g$ has the decay rate $t^{-(m+1)/2}$. This 
known observation gives information resembling our result, although it does not
give such a general decomposition of the initial data space as Theorem \ref{th1.0}. 
We also mention \cite{GW}, Appendix A, where analogous results in the form of 
spectral decompositions are derived  for more general equations. 
\ENR

There is an extensive literature dealing with  the decay rate of the solution to the Cauchy 
problem of the heat equation. For example, precise decay rates in the linear case
have been considered in \cite{C}, although most of the recent research is concentrated on
semilinear or other nonlinear generalizations of \eqref{1.2}---\eqref{1.3}. 
As a slightly random sample we mention the papers \cite{FWY,Fu,Gi,GV1,IK, Kaw,KST,PoQu,Qu,Zh}; 
see also  the monograph \cite{QS} for an
exposition. We especially mention the papers \cite{BKL,IK, IK1,IK2,T,T1}, where 
the asymptotic large time behavior of the semilinear problem is considered by 
separating the faster decay of terms with vanishing integrals. The paper  \cite{IK1}
contains the state of art in this direction and in fact has partially been a source
of inspiration for the present work.  

We organize our paper as follows. Section \ref{sec1.a} is devoted to the case $N=1$.
We discuss the known phenomenon that for special initial data with certain vanishing 
iterated integrals the decay rate can be made arbitrarily fast. This leads to the 
definition of special continuous linear functionals in the space $L_w^p(\bbR)$,
$1 \leq p < \infty$, and to the formulation
of Theorem \ref{th1.4} concerning the existence of Schauder basis annihilating
given functionals. We show how Theorem \ref{th1.0} follows from this result, although the
proof of  Theorem \ref{th1.4} is only presented in Section \ref{sec2}. Theorem \ref{th1.4} 
uses the concept of a shrinking Schauder basis: since an arbitrary  basis of the 
non-reflexive space  $L_w^1(\bbR)$ is not necessarily shrinking, this case requires a
separate treatment, which is contained in Lemma \ref{lem2.6}.

The case $N >1$
of  Theorem \ref{th1.0} will be considered in Section \ref{sec1.b}. The proof is based on 
decomposing a given initial data of several variables into a convergent sum of products of 
functions in one variable and using the already proven one-dimensional case. Here, our  
method requires the use of $L^1$-norms and a little abstract tensor product
techniques.  At the end of the Section \ref{sec1.b} we  discuss some interesting open 
problems.

We will use the following general notation. By $C,C'$ etc. we denote generic positive
constants, the exact value of which may change from place to place. The possible dependence,
say, on a parameter $p$ is indicated as $C_p$. By supp\,$f$ we denote the support of a function 
$f$ and by sp$(A)$ the linear span of a subset $A$ of a vector space. Its closure is denoted
by $\overline{{\rm sp}(A)}$. We write $\bbN  = \{ 1,2,
\ldots \}$, $\bbN_0 = \{0\} \cup\bbN$, and  $\bbR^\pm = \{ x \in \bbR \, : \, \pm x \geq 0 \}$. The characteristic or indicator function of a set $A$ is denoted by
$\mathbf{1}_A$. We use standard notation  
$L^p(\bbR^N)$, $L^p(0,1)$  etc. for unweighted Lebesgue spaces. 
Moreover, $X^*$ stands for  the dual of a Banach space $X$, i.e. the space
of bounded linear functionals on $X$. The norm of $X^*$ is denoted $\Vert \cdot \Vert_{X^*}$.
The identity operator $X \to X$ is denoted by ${\rm id}_X$.
For a linear operator $T$ between Banach spaces, $\Vert T \Vert $ denotes the operator norm. 

If $X$ denotes a separable Banach space over the scalar field $\bbK$ (either $\bbR$ or 
$\bbC$), we recall that a sequence $(e_n)_{n=1}^\infty 
\subset X$ 
is a Schauder basis, if every element $f \in X$ can be presented as a convergent
sum
\bea
f = \sum_{n=1}^\infty f_n e_n
\label{1.10}
\eea
where the numbers $f_n \in \bbK$ are unique for $f$. For example in a separable Hilbert space,
every orthonormal basis is a Schauder basis, but the converse is of course not true. 
There are many well-known constructions of Schauder bases in classical Banach spaces;
among them, the wavelet bases are most studied in the recent years. We refer to
\cite{LT}, \cite{Wo}  for this topic.

\section{Proof of Theorem \ref{th1.0} in the case $N=1$}
\label{sec1.a}

In this section we show how Theorem \ref{th1.0} follows from an abstract
result concerning bases which annihilate linear functionals, Theorem \ref{th1.4}. 
First, we recursively define the linear operators
\bea
& & I^0 f := f \ , \ \ I^m f(x) = \int\limits_{-\infty}^{x} I^{m-1} 
f(y) dy \ \ \mbox{for $m \in \bbN$} \ ,   \nonumber \\
& & J^m f = \int\limits_{-\infty}^0 
I^{m-1} f (y) dy \ \mbox{for 
$m \in \bbN$} . \label{1.16}
\eea
By the Cauchy formula for repeated integrations these can be written as 
\bea
I^m f(x)& = & \frac{1}{(m-1)!}\int\limits_{-\infty}^x (x-y)^{m-1} f(y) dy \ , 
\label{1.16m}  \\
J^m f & = &  \frac{1}{(m-1)!}\int\limits_{-\infty}^0 (-y)^{m-1} f(y) dy \ , \ \ m \in \bbN.
\label{1.16n}
\eea
The operators $I^m$ do not map $L_w^p(\bbR)$ even into $L^p(\bbR)$ (since 
$If:= I^1 f $ may be bounded from below by a positive constant for large $x$).
However, we have the following simple observations. We denote
$L_w^{p,-} (\bbR):= \{ f \in L_w^p(\bbR)  \, : \, {\rm supp} f \subset \bbR^-\} $.

\BEL
\label{lem1.2a}
$(i)$   If $m \in \bbN$ and $f \in L_w^p(\bbR)$, then the restriction of  $I^m f$ to 
$\bbR^{-}$ is  rapidly decreasing, as $x \to - \infty$: we have 
\bea
\sup\limits_{x \in \bbR^{-}} (1 + |x|)^k |I^m f(x)| 
\leq C_{k,m,p} \Vert f \Vert_{p,w} < \infty \ \ \forall \, k \in \bbN.
 \label{1.16c}
\eea

\noindent $(ii)$ If $m \in \bbN$ is given and $f \in L_w^{p,-}(\bbR)$ has the property that
$J^k f = 0$ for all $k\in \bbN$ with $k \leq m$, then 
\bea
{\rm supp}\, I^k f \subset \bbR^-
\label{1.16a}
\eea  
for all  $k \leq m$. In particular $I^k f \in L^1(\bbR)$ and
\bea
\Vert I^k f \Vert_1 \leq C_{k,p,w} \Vert f \Vert_ {p,w}  \label{1.16b}
\eea
for every $k \leq m$. 

\noindent $(iii)$  $J^m$ is a bounded linear functional 
on $L_w^{p}(\bbR)$. 
\ENL

Proof. As for  $(i)$, we consider  $p >1$ with the dual exponent $p'=p/(p-1)$. Then,
\eqref{1.16m}, \eqref{1.6}  and the H\"older inequality imply
for $x \leq 0$ 
\beas
 & & (1 + |x|)^k |I^m f(x)| \leq 
C_m (1 + |x|)^k 
\int\limits_{-\infty}^x |x-y|^{m-1} |f(y)| dy 
\rowleq
C_m (1 + |x|)^k 
\int\limits_{-\infty}^x |y|^{m-1}(1+ |y|)^{-k-m+1-\frac2{p'}} (1+ |y|)^{k+m-1+\frac2{p'}}|f(y)| dy 
\rowleq
C_m  \int\limits_{-\infty}^x (1+ |y|)^{-\frac2{p'}} (1+ |y|)^{k+m-1+\frac2{p'}}|f(y)| dy 
\rowleq
C_m  \Big(\int\limits_{-\infty}^0 (1+ |y|)^{-2} dy  \Big)^{1/{p'}}
\Big( \int\limits_{-\infty}^0  (1+ |y|)^{p(k+m-1+\frac2{p'})}|f(y)|^p dy \Big)^{1/p}
\rowleq
C_{k,m,p}
\Big( \int\limits_{-\infty}^0 w(y) |f(y)|^p dy \Big)^{1/p} \leq C_{k,m,p} \Vert f \Vert_{p,w}.
\eeas
The proof for the case $p=1$ is simpler, as the exponents $2/{p'}$ are omitted and the 
H\"older inequality is not needed. 

Concerning  $(ii)$, a  simple induction argument yields \eqref{1.16a}: assume that 
$J^k f = 0$ for 
all $k \leq m$ and  that $\tilde m < m$ and \eqref{1.16a} holds for all $k \leq \tilde m$. 
Then, by the definition of $I^{\tilde m + 1}$, for $x \geq 0$,
\beas
& &  I^{\tilde m + 1} f(x)
= 
\int\limits_{-\infty}^{0} I^{\tilde m} f (y) dy +
\int\limits_0^{x} I^{\tilde m} f(y) dy
\eeas
Here, the first term equals $J^{\tilde m + 1} f$ and is thus 0, and the second
term also vanishes by the induction assumption. 
The bound \eqref{1.16b} follows from \eqref{1.16a}, \eqref{1.16c} and an application
of the H\"older inequality. 


The statement $(iii)$ follows from \eqref{1.16n}, \eqref{1.6}, and  the H\"older inequality. \ \ $\Box$

\bigskip

The following fact about faster convergence rates for special initial data is known, but we 
need to present and prove a formulation, which  precisely fits to our arguments.

\BEP
\label{prop1.3}
Let $N=1$ and let $f \in L_w^{p,-}(\bbR)$ be such that for some $m \in \bbN$, 
it satisfies $J^k f =  0$ for all $k\in \bbN$ with $k \leq m$. Then, there holds the bound
\bea
\Vert e^{t \partial_x^2} f  \Vert_\infty \leq \frac{C_{p,w,m}\Vert f \Vert_{p,w}}{ t^{(1+m)/2}}   \label{1.17}
\eea
for the solution of \eqref{1.2}--\eqref{1.3} with the initial data $f$. 
\ENP

Proof. Let $m$ and $f$ be as in the assumption. We employ repeated integration by parts with 
respect to  $y$ in order to evaluate \eqref{1.4}. In this process there 
appear the expressions $I^kf$, which according to our assumptions and Lemma \ref{lem1.2a}
belong to $L^1(\bbR)$. 
At the first step we write
\beas
u(x,t) & =&  - \frac{1}{\sqrt{2 \pi t} }
\int\limits_{\bbR} \big( \partial_{y} e^{- \frac1{4t} (x-y)^2 } \big) (If)(y) dy 
= 
- \frac{1}{2\sqrt{2 \pi t}}
\int\limits_{\bbR} \frac{x-y}t e^{- \frac1{4t} (x-y)^2 }  (If)(y) dy ,
\eeas
where obviously the replacement term vanishes since both the Gaussian kernel and $If$ are 
rapidly  decreasing functions. 
Repeating integration by parts $k$ times, an induction proof shows that
\bea
u(x,t) & =&  
\frac{C}{\sqrt{t}}
\int\limits_{\bbR} P_k(x,y,t)  e^{- \frac1{4t} (x-y)^2 }  (I^kf)(y) dy 
\label{1.17b}
\eea
where the function $P_k$ is a finite sum of terms of the form
\bea 
C_{a}\frac{1}{t^{k/2}}\frac{(x-y)^{a}}{t^{a/2}}  \label{1.17c}
\eea
where $a \in \bbN_0$ and $C_a$ are some constants. Indeed, given $k \in \bbN$ and $a \geq 1$,
\beas
& & \partial_y \Big( \frac{1}{t^{k/2}}\frac{(x-y)^{a}}{t^{a/2}} 
e^{- \frac1{4t} (x-y)^2 } \Big)
\roweq
\Big(  \frac{-1}{t^{k/2}}\, \frac{a(x-y)^{a-1}}{t^{a/2}} + 
\frac{1}{t^{k/2}}\frac{2(x-y)^{a + 1}}{4 t^{a/2+1}} \Big) e^{- \frac1{4t} (x-y)^2 }
\roweq
\Big(  \frac{-1}{t^{(k+1)/2}}\, \frac{a(x-y)^{a-1}}{t^{(a-1)/2}} + 
\frac{1}{t^{(k+1)/2}}\frac{(x-y)^{a + 1}}{2 t^{(a+1)/2}} \Big) e^{- \frac1{4t} (x-y)^2 }
\eeas
where we have two terms of the form \eqref{1.17c} for $k+1$. 

We evaluate
\bea
& & \frac{1}{t^{1/2}}
\Big| \int\limits_{\bbR} \frac{1}{t^{k/2}}\frac{(x-y)^{a}}{t^{a/2}}   e^{- \frac1{4t} (x-y)^2 }  (I^kf)(y) dy  \Big|
\rowleq
\frac{C_{p,w,k,a}}{t^{(1+k)/2}} \sup\limits_{x,y\in \bbR }  \bigg( 
\Big( \frac{|x-y|^2}{4t} \Big)^{a/2}   e^{- \frac{|x-y|^2}{4t} }   \bigg)
\int\limits_{\bbR}
\big| (I^kf)(y) \big| dy 
\rowleq
\frac{C'_{p,w,k,a}}{t^{(1+k)/2}}  \Vert f \Vert_{p,w} ,
\eea
since the  supremum is bounded by a constant independent of $t$ and the integral is
bounded according to \eqref{1.16b}. We get the bound \eqref{1.17} for \eqref{1.17b}, by the 
remark on $P_k$ around \eqref{1.17c}. \ \ $\Box$

\bigskip

Given a Schauder basis $(e_n)_{n=1}^\infty$ of $X$ we denote for every $n \in \bbN$ by 
$P_n$ the basis projection
\bea
P_n f = P_n \Big(\sum_{k=1}^{\infty} f_k e_k \Big) = 
\sum_{k=1}^n f_k e_k \ , \ \mbox{ where}  \ 
f = \sum_{k=1}^{\infty} f_k e_k \in X.   \label{1.17m}
\eea
The number $K = \sup_n \Vert P_n \Vert$ is called  the basis constant of $(e_n)_{n=1}^\infty$;
the supremum defining $K$ is always finite, see \cite{LT}. 

\BED
\label{def1.5}
Let $x^* \in X^*$. We say that a Schauder basis $(e_n)_{n=1}^\infty$ of $X$  is {\em shrinking} with respect 
to $x^*$ if 
\bea
\lim_{n \to \infty} \Vert  x^* \circ ({\rm id}_X-P_n) \Vert_{X^*}  =0. \label{1.30}
\eea
\END

For a basis $(e_n)_{n=1}^{\infty}$ of $X$ consider the biorthogonal functionals $e_n^* \in 
X^*$, where $e_n^*(e_m) = \delta_{mn}$ (Kronecker delta); let $W = \overline{ {\rm sp} \{ 
e_n^* \, : \, n \in \bbN\} } \subset X^*$.  It is easily seen that  $(e_n^*)_{n=1}^{\infty}$ 
is a Schauder basis of $W$ with the basis projections $P_n^*$, where $P_n^*(x^*) = x^*\circ P_n$ for $x^* \in X^*$. However we have $W \not=X^*$ in general. 
We obtain that $(e_n)_{n=1}^{\infty}$ is shrinking with respect
to $x^* \in X^*$, if and only if $x^* \in W$.

Definition \ref{def1.5} extends slightly the classical notion of a shrinking basis, see  
\cite{LT}. A basis  $(e_n)_{n=1}^\infty$ of $X$ is shrinking, if it is
shrinking with respect to all elements in $X^*$ in the sense of the preceding definition, i.e. 
if $W=X^*$. In this case $X^*$ must be separable. It is well-known that every basis of $X$ 
is shrinking, if $X$ is reflexive. Again, see  \cite{LT} for more details.

Theorem \ref{th1.0} is a consequence of the following result, the proof of which is postponed 
to Section \ref{sec2}.

\BET
\label{th1.4}
Let $x_m^* \in X^*$ for all  $m \in \bbN$, and let  $\epsilon > 0$. 
Assume that $(\tilde e_n)_{n=1}^\infty$ is 	a Schauder basis of $X$ which is shrinking with respect to all 
$x^*_m$. Then, there exists an increasing  sequence $(n_m)_{m = 1}^\infty \subset \bbN$  and a basis $(e_n)_{n=1}^\infty$ of $X$ such that
\bea
x_m^*(e_n) = 0 \ \ \mbox{ for all } \ \  n \geq n_m.  \label{1.20}
\eea
If  $T:X \rightarrow X$ is the linear operator with $T \tilde e_n = e_n$ for all $n$, then we 
have
\bea
\Vert {\rm id}_X -T \Vert  < \epsilon.  \label{1.22}
\eea
\ENT

Obviously, condition \eqref{1.22} means that $T$ is a bijection and the new basis 
$(e_n)_{n=1}^\infty$ can be considered as perturbation of the given basis  $(\tilde e_n)_{n=1}^\infty$.

We repeat that every  Schauder basis of a Banach space $X$ is shrinking, if $X$ is reflexive. 
This is in particular true for any orthonormal basis in a Hilbert 
space. However, in order to treat the case $p=1$ we state the following result, which also 
will proven only in Section \ref{sec2}.

\BEL
\label{lem2.6}
There exists a Schauder basis $(\tilde{e}_n^-)_{n=1}^{\infty}$ of
$L_{w}^{1,-}(\mathbb{R})$ which is shrinking for all  functionals $J^m$ defined in (2.3).
\ENL

Proof of Theorem \ref{th1.0}. Let $p$ and $w$ be as in the assumption and first consider 
the  Banach space $X= L_w^{p,-}(\bbR)$. The functionals $J^m =: x_m^* $ of \eqref{1.16} are 
well defined and bounded on $X$, by Lemma \ref{lem1.2a}, $(iii)$. We fix a basis
$(\tilde e_n^-)_{n=1}^\infty$, which is shrinking with respect to all $x_m^*$; in the
case $p=1$ we use Lemma \ref{lem2.6} to find this. Then, Theorem \ref{th1.4} yields the 
desired basis $(e_n^-)_{n=1}^\infty$ of $L_w^{p,-}(\bbR)$ and  the sequence of indices $(n_m)_{m=1}^\infty$; in 
particular,  given  $m \in \bbN$  we have 
\bea
J^k(e_n^-) = 0  \label{1.24}
\eea
for every $k \leq m$, $n \geq n_m$.  To see that \eqref{1.12a} holds for a given $m$ and for any 
initial data  $f^- \in G_{n_m}^- := \overline{ {\rm sp}\, \{ e_n^- \, : \, n \geq n_m \} } \subset L_w^{p,-}(\bbR) $  we remark that such a $f^-$ has a representation
\bea
f^-=\sum_{n = n_m}^\infty f_n^- e_n^- .
\eea
Since this series converges in $L_w^p(\bbR)$ and every $J^k$ is a continuous mapping, \eqref{1.24} 
implies $ J^k g = 0 $ for all $ k \leq m$. Hence, \eqref{1.12a} follows from Proposition
\ref{prop1.3}. 

To complete the proof we remark that the space $L_w^p(\bbR)$ equals in a natural way the 
direct sum  $L_w^{p,-} (\bbR) \oplus L_w^{p,+}(\bbR)$, where the second component is defined as the closed 
subspace of $L_w^p(\bbR)$ consisting of functions with supports in 
$\bbR^{+}$. The functions 
\bea
e_n^+ := e_n^- \circ \psi 
\ \ , \ \ \mbox{where $\psi(x) := -x 
\ \forall x\in \bbR$}
\eea
form a Schauder basis of $ L_w^{p,+}(\bbR)$, which plays the same role as the basis 
$(e_n^-)_{n=1}^\infty$
has in $L_w^{p,-}(\bbR)$. This follows from the formal commutation relations
\bea
\partial_x^2 ( f \circ \psi ) = (\partial_x^2 f ) \circ\psi \ \ , \ \ 
e^{t \partial_x^2 } (f \circ \psi) = \big( e^{t \partial_x^2 } f \big)  \circ \psi. 
\eea

Consequently, the union of the sequences $(e_n^-)_{n=1}^\infty$ and $(e_n^+)_{n=1}^\infty$ is 
the desired Schauder basis.  \ \ $\Box$

\section{Proofs of Theorem \ref{th1.4} and Lemma \ref{lem2.6}.}
\label{sec2}

We need the following elementary
\BEL
\label{lem2.1}
Let $(h_n)_{n=1}^\infty$ be a basis of the Banach space $Y$ with basis projections $Q_n$, $n \in \bbN$, and 
basis constant $K$. Moreover, let $T:Y \rightarrow Y$ be a linear operator with 
$c:= \Vert {\rm id}_Y -T \Vert < 1$. Then $(Th_n)_{n=1}^\infty$ is a basis of $Y$ with basis constant at 
most $K(1+c)/(1-c)$.
\ENL

Proof. By the assumption and the Neumann series,  $T$ is an isomorphism (linear homeomorphism),
and we have $T^{-1} = \sum_{k=0}^{\infty} ({\rm id}_Y-T)^k$, hence $\Vert T^{-1}\Vert 
\leq (1-c)^{-1}$. Moreover, $\Vert T\Vert \leq 1+ \Vert  {\rm id}_Y -T\Vert   \leq 1+c$. 
Hence, $(Te_n)_{n=1}^\infty$ is a basis of $Y$ with basis projections $T P_n T^{-1}$ and basis constant
at most  $ K \Vert  T\Vert \,  \Vert  T^{-1}\Vert   \leq K(1+c)/(1-c)$. \ \ $\Box$

\BEP
\label{prop2.2}
Let $(\tilde h_n)_{n=1}^\infty$ be a basis of the Banach space $Y$ with basis projections 
$Q_n$, $n \in \bbN$, and 
basis constant $K$. Moreover, let $L,M \in \bbN$ and assume that $y^*_m \in Y^*$, $m \in 
\bbN$, satisfy
\beas 
y_1^*|_{({\rm id}_Y-Q_L)Y}, \ldots, y_M^*|_{({\rm id}_Y-Q_L)Y} = 0
\eeas
and
\bea 
\lim_{n \rightarrow \infty} \Vert  y^*_m|_{({\rm id}_Y-Q_n)Y}\Vert   = 0 \ \ \mbox{ for all } m. 
\label{2.1} 
\eea
Then for any $\delta > 0$ there is a basis $(h_n)_{n=1}^\infty$ of $Y$ and an index $N > L$ with
\bea
& &  h_n = \tilde h_n , \ \ n=1, \ldots,N, \nonumber \\
& & y^*_{M+1}( h_n)=0 \ \ \mbox{ if } n > N,  \ \ \ y^*_k( h_n)=0 \ \ \mbox{ for } 
k=1, \ldots, M \ \mbox{and} \ n \geq L+1,
\label{2.2}
\eea
and 
\bea
\Vert  {\rm id}_Y -S\Vert    \leq K \delta
\label{2.3}
\eea
for the linear operator $ S:Y \rightarrow Y$  with  $S \tilde h_n= h_n$ for all 
$n \in \bbN$. 
The basis constant of $(h_n)_{n=1}^\infty$ is at most $K(1+K\delta)/(1-K\delta)$.
\ENP

Proof. If $y^*_{M+1}|_{({\rm id}_Y-Q_L)Y} =0$ then we can take $h_n = \tilde h_n$ for all $n$.
Otherwise let $N > L$ be large enough and put
\beas
\rho = \frac{\Vert  y^*_{M+1}|_{({\rm id}_Y-Q_N)Y}\Vert  }{\Vert  y^*_{M+1}|_{(Q_N -Q_L)Y}\Vert  }. 
\eeas
According to \eqref{2.1} we can choose $N$ so large that $\rho < \delta$ and 
\bea
K \rho < 1. \label{2.4}
\eea
In fact $\rho$ can be made arbitrarily small since the denominator in 
the definition of $\rho$ goes to $\Vert y^{*}_{M+1}|({\rm id}_Y-Q_L)Y \Vert>0$ if $N$ 
tends to $\infty$ while the numerator tends to 0 in view of \eqref{2.1}.
We find $x \in (Q_N- Q_L)Y$ with $\Vert  x\Vert   =1$ and $y^*_{M+1}(x) = 
\Vert   y^*_{M+1}|_{(Q_N-Q_L)Y}\Vert  $. (Take into account that $(Q_N-Q_L)Y$ is finite 
dimensional.) 
	 	
Put $Sf= f$ if $f \in Q_NY$ and 
\bea
Sg = g - \frac{y^*_{M+1}(g)}{ \Vert   y^*_{M+1}|_{(Q_N-Q_L)Y}\Vert  } x \ \ \ \mbox{ if } g \in 
({\rm id}_Y-Q_N)Y.
\label{2.5}
\eea
Then we have
\bea
\Vert  f+g - S(f+g)\Vert   = \Vert  g-Sg\Vert   \leq \rho \Vert  g\Vert   
\leq \rho K \Vert  f+g\Vert  . 
\label{2.6}
\eea
Let $h_n = S\tilde h_n$ for all $n$. According to Lemma \ref{lem2.1} and in view of \eqref{2.4}, 
$(h_n)_{n=1}^\infty$ is a basis of $Y$ with  basis constant smaller than or equal to
\beas
K \left(\frac{1+K \rho}{1-K\rho} \right) \leq K \left(\frac{1+K \delta}{1-K\delta} \right).
\eeas
Formula \eqref{2.5} yields $y^*_{M+1}(h_j) = 0 $ if $j > N$. Moreover, since 
$x\in ({\rm id}_Y-Q_L)Y$ we have $y^*_k(h_l)=0$ for $k \leq M$, $l \geq L+1$.
Together with \eqref{2.6} this proves the proposition. $\Box$

\bigskip

Conclusion of the proof of Theorem \ref{th1.4}. Consider $\delta_n > 0$ such that
\beas
\sum_{n=1}^{\infty}K 2^{n-1}\delta_n \leq \epsilon, \ \ \frac{1+K2^{n-1} \delta_n}{1-K2^{n-1}\delta_n} \leq 2 \ \ \mbox{ and } \ \ \prod_{n=1}^{\infty}
\left( \frac{1+K2^{n-1} \delta_n}{1-K2^{n-1}\delta_n} \right) \ \ \mbox{ converges. }
\eeas
Then, we use induction and apply Proposition \ref{prop2.2} as follows. 

We start with the basis $(\tilde e_n)_{n=1}^\infty =: (e_n^{(1)})_{n=1}^\infty$ and
$n_1 := 0$. 
If we are in the step $m$, and we already have the indices $n_k$, $k \leq m$, 
and a basis  $(e_n^{(m)})_{n=1}^\infty$ with basis constant at most
\beas 
K \prod_{k=1}^{m-1}
\left( \frac{1+K2^{k-1} \delta_k}{1-K2^{k-1}\delta_k} \right) \leq 2^{m-1} K ,
\eeas
such that 
$ 
x_k^* (e_n^{(m)} )= 0
$ 
for all $n \geq n_k$ and all $k \leq m$, then we apply Proposition \ref{prop2.2} with 
$\tilde h_n = e_n^{(m)}$, $L= n_m$, $M=m$ and
$\delta =\delta_m$. This yields  an index $N >  n_m$ and a basis 
$(e_n^{(m+1)})_{n=1}^\infty$
with basis constant not larger than
\beas 
K \prod_{k=1}^{m} \left( \frac{1+K2^{k-1} \delta_k}{1-K2^{k-1}\delta_k} \right)
\eeas
such that, in view of \eqref{2.2}, $e_n^{(m+1)} = e_n^{(m)}$ for $n \leq N$ and 
$ 
x^* (e_n^{(m+1)} )= 0
$ 
for all $n > N $. Put $n_{m+1} = N$ and continue the induction. 
	  
At the  $m$th step of the process, the first $n_m$ elements of the basis remain unchanged  so that we end up with a basic sequence $(e_n)_{n=1}^\infty$ with basis constant
at most 
\beas 
K \prod_{k=1}^{\infty} \left( \frac{1+K2^{k-1} \delta_k}{1-K2^{k-1}\delta_k} \right)  
\eeas
and such that \eqref{1.20} holds. 
In view of  \eqref{2.3} the linear operator $T: X \rightarrow X$ with $T \tilde e_n= e_n$ for all $n$ 
satisfies
\beas
\Vert {\rm id}_Y-T\Vert   \leq \sum_{m=1}^{\infty}K
\prod_{k=1}^{m-1} \left( \frac{1+K2^{k-1} \delta_k}{1-K2^{k-1}\delta_k} \right) \delta_n 
\leq
\sum_{m=1}^{\infty}K2^{m-1}\delta_m \leq \epsilon. 
\eeas
If we choose $\epsilon < 1$ then $T$ is surjective and  $({e}_n)_{n=1}^\infty$ is a basis of $X$
with the required properties.	\ \   $\Box$

\bigskip

Proof of Lemma \ref{lem2.6}.
We consider the Haar system $(e_n)_{n=1}^{\infty}$ in $L^1(0,1)$, where  $e_1 \equiv 1$
and 
$$ e_{2^k+j}(t)= \left\{ \begin{array}{rl}
1  , & \ \mbox{ if } t \in [(2j-2)2^{-k-1},(2j-1)2^{-k-1}], \\
-1 , & \ \mbox{ if }t \in [(2j-1)2^{-k-1},(2j)2^{-k-1}], \\
0  , & \ \mbox{ otherwise, }   
\end{array}
\right. 
$$ 
for $k = 0,1,2, \ldots$ and $j=1, \ldots, 2^k$.
 It is well-known that the Haar system is a Schauder basis for $L^1(0,1)$ with basis constant 1 (see \cite{LT}).
Put 
$$ 
A_{2^k+j-1}= \Big[ \frac{j-1}{2^{k}}, \frac{j}{2^{k}} \Big]. 
$$
Then we have $\mathbf{1}_{A_1}= e_1$, $\mathbf{1}_{A_2}= (e_1+e_2)/2$ and 
$\mathbf{1}_{A_3}=(e_1-e_2)/2$. By induction we see that any element $\mathbf{1}_{A_m}$ 
is a linear combination of the Haar functions $e_n$.

For $h \in L^{\infty}(0,1)$ let $\Phi_h$ be the linear functional on $L^1(0,1)$ 
defined by 
$$ 
\Phi_h(f) = \int\limits_0^1 f(s) h(s) ds \ \ \ \mbox{ for all }  f \in L^1(0,1). 
$$
Recall that the map $h \mapsto \Phi_h$ is an isometric isomorphism between $  L^{\infty}(0,1)$ 
and $L^1(0,1)^*$. It is easily seen that the biorthogonal functionals $e_n^*$ of the 
Haar elements $e_n$ are, up to  constant factors,  the functionals $\Phi_{e_n}$. Let 
$W = \overline{ {\rm sp}\, \{ \Phi_{e_n} \} } \subset L^1(0,1)^*$.
Then $\Phi_h \in W$ for any linear combination $h$ of the functions $\mathbf{1}_{A_n}$.

Define $\alpha : \  ]0,1] \rightarrow \ ] - \infty,0]$ by
$\alpha(s) = \log s$, $s \in ]0,1]$, and
$(Sf)(s)= f( \alpha(s))w(\alpha(s))/s$ for
$ f \in L^{1,-}_{w}(\mathbb{R})$. Then $S$ is
an isometric isomorphism between $L^{1,-}_{w}(\mathbb{R})$ and $L^1(0,1)$. In particular we 
have
$$ 
J^mf = \int\limits_0^1 (Sf)(s) \frac{(-\alpha(s))^{m-1}}{(m-1)!w(\alpha(s))} ds
\ \ \ \mbox{ for all }m \in \bbN. 
$$
With 
$$ 
g(s) = \frac{(-\alpha(s))^{m-1}}{(m-1)!w(\alpha(s))}, \ \ s\in ]0,1], \ \ \mbox{ and } \ \ g(0)=0 
$$
we obtain 
\bea
J^mf = \Phi_g(Sf).  \label{2.31}
\eea
In view of (1.5) the function $g$ is continuous and hence 
uniformly continuous  on $[0,1]$. This means that, for any $ \epsilon > 0$,  we find $k$ 
and a linear combination $g_{\epsilon }$ of  the characteristic functions 
$\mathbf{1}_{A_{2^k+j-1}}$, $j=1, \ldots, 2^k$, such that
$$ 
\Vert g - g_{\epsilon }\Vert_{\infty} \leq \epsilon,  
$$
and hence $\Vert \Phi_g - \Phi_{g_{\epsilon }} \Vert  \leq \epsilon $. Since $\Phi_{g_{\epsilon }}\in W$ for all $\epsilon$
we conclude $\Phi_g \in W$.

Finally, let $ \tilde{e}_n^- = S^{-1}e_n$ for all $n$. Then  $(\tilde{e}_n^-)_{n=1}^{\infty}$ is a Schauder basis of
$L_{w}^{1,-}(\mathbb{R})$. The norm-closed linear span of the biorthogonal functionals is equal to
$S^*W = \{ w^*\circ S : w^* \in W \}$. From \eqref{2.31} we obtain $J^m \in S^*W$ for all $m \geq 1$ and therefore $(\tilde{e}_n^-)_{n=1}^{\infty}$ is shrinking for the
functionals  $J^m$ (see the remark after Definition \ref{def1.5}). $\Box$

\section{Proof of Theorem \ref{th1.0} in the case $N > 1$}
\label{sec1.b}

Returning to the proof of Theorem \ref{th1.0}, when $N> 1$, we assume that the weight 
$w : \bbR^N \to \bbR^+$ and $m \in \bbN$ are given. 
First we select a weight $\tilde w : \bbR^N \to \bbR^+$ such that 
$ L_w^p(\bbR^N) \subset L_{\tilde w}^1(\bbR^N)$ and such that 
\bea
\tilde w (x) = \prod_{j=1}^N v(x_j) \label{5.01}
\eea
for some  weights $v$ on $ \bbR$  satisfying the assumptions around \eqref{1.6}
in the one-dimensional case. One can  for example find $\tilde w$  as follows. Define for 
every $k \in \bbN$ the number $B_k$ such that
\bea
B_k = \sup\limits_{x \in \bbR^N} \frac{(1+|x|)^k}{w(x)}
\eea
and then set
\bea
\tilde v(x) = \prod_{j=1}^N \sum_{k=1}^\infty 2^{-k} B_k^{-1/N} (1 + |x_j|)^{k/N}.
 \label{5.02}
\eea
We then have
\bea
\tilde v(x) & \leq &  \prod_{j=1}^N \sum_{k=1}^\infty 2^{-k} B_k^{-1/N} (1 + |x|)^{k/N}
\rowleq
\prod_{j=1}^N \sum_{k=1}^\infty 2^{-k} (1 + |x|)^{k/N}
\inf\limits_{y \in \bbR^N} \frac{w(y)^{1/N}}{(1+|y|)^{k/N}}
\leq 
\prod_{j=1}^N \sum_{k=1}^\infty 2^{-k}w(x)^{1/N} =  w(x)
\eea
If  $p=1$, we take $\tilde w = \tilde v$, and for $p >1$ we set
\bea
\tilde w(x) =  \tilde v(x)^{1/p}  \prod_{j=1}^N (1 + |x_j|)^{-2/p' }  \label{5.4}
\eea
where  $p' = p/(p-1)$ is the dual exponent of $p$.
Then, $\tilde w$ is still as in \eqref{1.6}, and moreover, for every $f \in 
L^p_{\tilde v} (\bbR^N)$ we have by the H\"older inequality 
\bea
& & \Vert f \Vert_{1,\tilde w} = \int\limits_{\bbR^N} |f| \tilde w dx \leq 
\int\limits_{\bbR^N} |f(x)| \tilde v(x)^{1/p}  \prod_{j=1}^N (1 + |x_j|)^{-2/p' }  dx
\rowleq 
\Big( \int\limits_{\bbR^N} |f(x)|^p  \tilde v(x) dx \Big)^{1/p}
\Big( \int\limits_{\bbR^N}  \prod_{j=1}^N (1 + |x_j|)^{-2 }  dx \Big)^{1/p'}
\leq C \Vert f \Vert_{p,\tilde v}
\eea
so that  $ L^p_{w}(\bbR^N) \subset L^p_{\tilde v}(\bbR^N) 
\subset L^1_{\tilde w}(\bbR^N)$, and $\tilde w$ is of the form \eqref{5.01} 
with 
\beas
v(x_j) = (1 + |x_j|)^{-2/p'} \sum_{k=1}^\infty 2^{-k} B_k^{-1/N} (1 + |x_j|)^{k/N} .
\eeas

Using Theorem \ref{th1.0} with the weight $v$ (in the place  of 
$w$) we find the Schauder basis $(\tilde e_n)_{n=1}^\infty$ in $L_{v}^1(\bbR)$
and the increasing sequence  $(n_m)_{m=1}^\infty$ such that
\eqref{1.12a} holds. We may and do require that the basis $(\tilde e_n)_{n=1}^\infty$ is 
normalized in $L_{v}^1(\bbR)$ so that  $\Vert \tilde e_n \Vert_{1, v} = 1$ for every $n$.
We denote the corresponding $n$th basis projection \eqref{1.17m} by $\tilde P_n$
and the corresponding complementary projection $\tilde Q_n = {\rm id}_{L_v^1(\bbR)} -
\tilde P_n$.

Since the weighted Lebesgue measure $\tilde w  dx$ on $\bbR^N$ is the product of the 
$N$ measures $v dx_j$ on $\bbR$ by \eqref{5.01}, we can apply the theory 
tensor product norms and present the space $L^1_{\tilde w}(\bbR^N)$ as the $N$-fold 
projective tensor product
\bea
L^1_{\tilde w}(\bbR^N) = L^1_v(\bbR) \widehat \otimes_\pi \ldots  \widehat \otimes_\pi
L^1_v(\bbR) ,  \label{5.6pp}
\eea
see \cite{Tr}, Ch. 46 and in particular Exercise 46.5. We need a few facts concerning 
\eqref{5.6pp}: according to the  definition of the projective tensor product, every 
$f \in L^1_{\tilde w}(\bbR^N)$ can be written as 
\bea
f(x) = \sum_{k=1}^\infty \lambda_k f^{(1,k)}(x_1)  f^{(2,k)}(x_2) \ldots  f^{(N,k)} (x_N)
\label{5.6}
\eea
where  the numbers  $\lambda_k$ form an absolutely summable sequence,
\bea
\sum_{k=1}^\infty |\lambda_k| \leq C \Vert f \Vert_{1,\tilde w}  \label{5.8}
\eea
and every $f^{(j,k)}$ belongs to the space $ L^1_{v}(\bbR)$ and has the bound
\bea
\Vert f^{(j,k)} \Vert_{1, v} \leq 1 .  \label{5.10}
\eea
It follows that the sum \eqref{5.6} converges absolutely in the space
$L_{\tilde w}^1(\bbR^N)$.

The second fact is that the functions 
\bea
e_{\bar n} = \tilde e_{n(1)} \otimes \ldots \otimes \tilde e_{n(N)}
\ \ (\mbox{i.e.} \ e_{\bar n} (x) = \tilde e_{n(1)} (x_1) \ldots \tilde e_{n(N)} (x_N), 
\ x \in \bbR^N ) \label{5.10y}
\eea
where $\bar n = \big( n(1), n(2), \ldots , n(N)\big) \in 
\bbN^N$ runs over all $N$-tuples, form a Schauder basis of the space 
$ L^1_{\tilde w}(\bbR^N)$, see for example \cite{Ho}. We will prove the theorem by using 
this basis. Here, we  do not need to order
the basis explicitly with an index in $\bbN$; nevertheless,  every operator 
\beas
P_n := \tilde P_n \otimes \ldots \otimes \tilde P_n \ \ , \ \ n \in \bbN ,
\eeas
is a basis projection with $\sup_n \Vert P_n \Vert < \infty $ in the operator norm
of $ L^1_{\tilde w}(\bbR^N)$, although not all basis projections of the basis
$(e_n)_{n=1}^\infty$ are of this form. The complementary projection can be 
written as a finite sum
\bea
Q_n := I - P_n = 
\sum_{\sigma  \in \bbS} \tilde R_{n,\sigma(1)} \otimes \ldots \otimes \tilde
R_{n,\sigma(N)} \ \ , \ \ n \in \bbN \label{5.13b}
\eea
where we write $I$ for the identity operator on $ L^1_{\tilde w}(\bbR^N)$ for brevity,
\bea
\sigma = ( \sigma(1), \ldots , \sigma(N) ) \in \bbS := \prod_{j=1}^N \{1,2\}
\smallsetminus \{(1,1, \ldots ,1)\}
\eea
and 
\begin{equation}
\tilde R_{n,k} = \left\{
\begin{array}{ll}
\tilde P_n , & \ \ \ \ k=1 \\
\tilde Q_n , & \ \ \ \  k=2 .
\end{array}
\right.
\end{equation} 
In other words, the sum \eqref{5.13b} consists of exactly those terms, where at least one factor equals $\tilde Q_n$; the sum has $2^N -1$ terms. 

\BEL
\label{lem5.9}
If $f \in  L^1_{\tilde w}(\bbR^N)$ and $P_n f = 0$ for some $n \in \bbN$, then 
$f$ has a representation \eqref{5.6}, where for every $k$ at least one factor
$f^{(j,k)}$ satisfies 
\bea
\tilde P_n f^{(j,k)} = 0 \ \ , \ \ \mbox{equivalently}, \ \ 
f^{(j,k)} =  \tilde Q_n f^{(j,k)} ,  \label{5.13c}
\eea
and  the bounds \eqref{5.8} and \eqref{5.10} still hold true. 
\ENL

Proof. If $f$ is given and  \eqref{5.13c} does not already hold for its representation
\eqref{5.6}, we write using the definition of the tensor product operator
\eqref{5.13b} and the absolute convergence of the series \eqref{5.6}
\bea
f(x) & = & Q_n  f(x) = \sum_{k=1}^\infty \lambda_k 
\sum_{\sigma \in \bbS} \tilde R_{n, \sigma(1)} f^{(1,k)}(x_1) 
\ldots   \tilde R_{n, \sigma(N)} f^{(N,k)}  (x_N)
\roweq
\sum_{k=1}^\infty \sum_{\sigma \in \bbS}
B^N \lambda_k 
\frac1B \tilde R_{n, \sigma(1)} f^{(1,k)}(x_1)  
\ldots  \frac1B \tilde R_{n, \sigma(N)} f^{(N,k)} (x_N)  ,
\label{5.13d}
\eea
where $B >0 $ is a uniform bound for the operator norms of the projections 
$\tilde P_n$ and $\tilde Q_n$ in the space $L_v^1 (\bbR)$, and the double sequence
$\big( B^N \lambda_k \big)_{k\in \bbN, \sigma \in \bbS} $ is still absolutely
summable,  since $\bbS$ has the fixed  number $2^{N-1}$ of  terms.
Thus \eqref{5.13d} is the desired representation of $f$.   \ \ $\Box$

\bigskip

We show that the basis \eqref{5.10y} satisfies the claim of Theorem \ref{th1.0}.
Let $m \in \bbN$ be given and let $n_m$ be as chosen above; we write for brevity
$n_m=:n$. As remarked above, $P_n$ is a basis projection related to the
basis \eqref{5.10y}. We assume that $f \in L_{\tilde w}^1 (\bbR^N)$ is such that 
$P_{n} f =0$, and take a representation \eqref{5.6} with the properties given
by Lemma \ref{lem5.9}. We consider an arbitrary  term of \eqref{5.6}, with some
abuse of notation in the variables:
\beas
& & e^{t \Delta} \bigg( \frac1B \tilde R_{n, \sigma(1)} f^{(1,k)}  (x_1)
\ldots  \frac1B \tilde R_{n, \sigma(N)} f^{(N,k)} (x_N)  \bigg)
\roweq
B^{-N} e^{t \Delta} \bigg(
\Big(  \prod_{j:\sigma(j) = 1} \tilde R_{n, \sigma(j)}  f^{(j,k)}(x_j)
\prod_{j:\sigma(j) = 2} \tilde R_{n, \sigma(j)}  f^{(j,k)} (x_j)\Big)  \bigg)
\roweq
B^{-N}  \prod_{j:\sigma(j) = 1} e^{t \partial_{x_j}^2}\tilde R_{n, \sigma(j)}  f^{(j,k)}
(x_j) 
\prod_{j:\sigma(j) = 2}  e^{t \partial_{x_j}^2} \tilde R_{n, \sigma(j)}  f^{(j,k)} (x_j)
\eeas
Here, if $\sigma(j)=1$, we have $ \tilde R_{n, \sigma(j)} = \tilde P_n$, and thus
by the uniform boundedness of the operator norms of $\tilde P_n$ in $L_v^1(\bbR)$ and 
\eqref{5.10},
\bea
\big\Vert e^{t \partial_{x_j}^2}\tilde R_{n, \sigma(j)}  f^{(j,k)} \big\Vert_{\infty} 
\leq \frac{C}{t^{1/2}}\big\Vert \tilde R_{n, \sigma(j)}  f^{(j,k)}
\big\Vert_{1,v} \leq  \frac{C'}{t^{1/2}} . 
\eea
However, if $\sigma(j)=2$, we have $ \tilde R_{n, \sigma(j)} = \tilde Q_n$,
and thus the choices made above and \eqref{1.12a} imply
\bea
& & \big\Vert e^{t \partial_{x_j}^2}\tilde R_{n, \sigma(j)}  f^{(j,k)} \big\Vert_{\infty} 
= \big\Vert e^{t \partial_{x_j}^2}\tilde Q_n  f^{(j,k)} \big\Vert_{\infty} 
\rowleq 
\frac{C}{t^{m} } \big\Vert \tilde Q_n   f^{(j,k)}
\big\Vert_{1,v} \leq \frac{C'}{t^{m} }
\eea
for $t\geq 1$.

By Lemma \ref{lem5.9}, every term in \eqref{5.6} has at least one factor with
$\sigma(j) =2$, hence, 
\beas
 & & \big\Vert e^{t \Delta } f \big\Vert_\infty
\leq 
\sum_{k=1}^\infty |\lambda_k| B^N 
\bigg\Vert e^{t \Delta} \bigg(\frac1B \tilde R_{n, \sigma(1)} f^{(1,k)}(x_1)
\ldots   \frac1B \tilde R_{n, \sigma(N)} f^{(N,k)} (x_N)  \bigg) \bigg\Vert_\infty
\rowleq
C \sum_{k=1}^\infty |\lambda_k|  \Vert  f \big\Vert_{1,\tilde w} 
\frac{1}{t^{m+(N-1)/2} }
\leq  \frac{C'}{t^{m +(N-1)/2} } \Vert  f \big\Vert_{1,\tilde w} 
\leq  \frac{C'}{t^{m} } \Vert  f \big\Vert_{1,\tilde w} 
\eeas
for $t \geq 1$. \ \ $\Box$

\bigskip

We conclude by a discussion. First, we remark that in the case $p=2$ it is possible to use 
standard Hilbert space methods (Fr\'echet-Riesz theorem and Gram-Schmidt method) and 
give an existence proof for an orthonormal basis in $L_w^{2, \pm}(\bbR)$ with the 
property \eqref{1.20}  for the functionals $x_m^* := J^m$. 
This yields the existence of an {\it orthonormal} 
basis in Theorem \ref{th1.0}.$1^\circ$, for $L_w^{2}(\bbR)$.


The heat equation is a classical albeit simplified model for the heat conduction or
linear diffusion processes. 
Since our discovery  is basically an existence proof, its  possible physical relevance 
depends on concrete examples of Schauder basis and estimates of the magnitude of the 
numbers $n_m$. We pose the problem:

\smallskip

\noindent $1^\circ$. {\it Given a weight $w$, $1\leq p < \infty$ and $m \in \bbN$, minimize 
the  number $n_m$ in Theorem \ref{th1.0}.}

\smallskip

\noindent Of course, the result in higher dimensions should be improved. 

\smallskip

\noindent $2^\circ$. {\it Find a Schauder basis in the space $L_w^p(\bbR^N)$, 
$N >1$, with the same properties as in  Theorem \ref{th1.0}, $1^\circ$.}

\smallskip

\noindent Finally, we ask if it is possible in the case $N >1$ to use the one-dimensional 
result in such a way that  the weight only needs to be fast growing in one coordinate
direction and milder assumptions are sufficient in other directions. For example:

\smallskip

\noindent $3^\circ$. {\it Does Theorem \ref{th1.0}, $2^\circ$,  hold for the weight
\bea
w(x) = e^{|x_1|} \prod_{j=2}^N (1 + |x_j|)^2 .
\eea
}

\bigskip

{\bf Acknowledgements.} The authors would like to thank Thierry Gallay (Grenoble) for 
discussions which helped in   the final formulation of our results. 
The research of Bonet was partially supported by the projects MTM2016-76647-P and GV Prometeo 2017/102. The research of Taskinen was partially supported by the research
grant from the Faculty of Science of the University of Helsinki.

\bigskip

\noindent \textbf{Authors' addresses:}%
\vspace{\baselineskip}%

\noindent
Jos\'e Bonet: Instituto Universitario de Matem\'{a}tica Pura y Aplicada IUMPA,
Universitat Polit\`{e}cnica de Val\`{e}ncia,  E-46071 Valencia, Spain

\noindent
email: jbonet@mat.upv.es \\

\noindent
Wolfgang Lusky: FB 17 Mathematik und Informatik, Universit\"at Paderborn, D-33098 Paderborn, Germany.

\noindent
email: lusky@uni-paderborn.de \\

\noindent
Jari Taskinen: Department of Mathematics and Statistics, P.O. Box 68,
University of Helsinki, 00014 Helsinki, Finland.

\noindent
email: jari.taskinen@helsinki.fi

\end{document}